%% file: ClasseDeConjugaison.tex
\documentclass[dvips]{article}
\usepackage{amsmath,amssymb,amsfonts,amsthm,epsfig,latexsym,makeidx,graphics}
\usepackage[usenames,dvipsnames]{color}
 \def\NN{{\mathbb N}}  
\def\QQ{{\mathbb Q}} \def\RR{{\mathbb R}}  
   
 \def\ZZ{{\mathbb Z}}

\def\De{\Delta}

   \def\cN{{\cal N}} 
\def\cC{{\cal C}}

\newtheorem{theo}{Theorem}[section]

\newtheorem{clai}{Claim}

\newtheorem{lemm}{Lemma}[section]
\newtheorem{coro}[lemm]{Corollary}

\newtheorem{defi}[lemm]{Definition}
\newtheorem{prop}[lemm]{Proposition}

\newtheorem{ques}{Question}

\newtheorem{rema}[lemm]{Remark}

\newenvironment{demo}[1][Proof]{\noindent {\bf #1~: }}{\hfill$\Box$\medskip}

\title{Smooth Conjugacy classes of circle diffeomorphisms with irrational rotation number}
\author{C. Bonatti and N. Guelman}

\begin{document}
\maketitle

\renewcommand{\thefootnote}{}



\renewcommand{\thefootnote}{\arabic{footnote}}
\setcounter{footnote}{0}%
%

\begin{abstract} In this paper we prove the $C^1$-density of every $C^r$-conjugacy class in the closed subset of diffeomorphisms of the circle with a given irrational rotation number.
\end{abstract}

\input{intro}

\input{Geometrie}

\input{Denjoy}

\input{Prova}
\input{Characteristic}
\input{Perturbations}

\input{biblio}
\noindent Christian Bonatti,

\noindent {\small Institut de Math\'ematiques de Bourgogne\\
UMR 5584 du CNRS}

\noindent {\small Universit\'e de Bourgogne, Dijon 21004, FRANCE}

\noindent {\footnotesize{E-mail : bonatti@u-bourgogne.fr}}

\vskip 5mm

\noindent Nancy Guelman

\noindent {\small I.M.E.R.L. Facultad de Ingenier\'ia.\\  Universidad de la Rep\'ublica}

\noindent {\small C.C. 30, Montevideo, Uruguay.}

\noindent{\footnotesize{E-mail:  nguelman@fing.edu.uy}}

\end{document}

%% file: intro.tex
\section{Introduction}

One knows from H. Poincar\'e that the dynamic of a homeomorphism $f$ of the circle depends strongly on the rotation number $\rho(f)$: the existence of periodic orbits is equivalent  to the rationality of $\rho(f)$.  If, on the contrary, the rotation number is irrational then $f$ is semi-conjugated to the corresponding irrational rotation. The non-injectivity of the semiconjugacy consists in colapsing each wandering interval to a point. In the thirties, A. Denjoy exhibited examples of $C^1$-diffeomorphisms with irrational rotation number but having wandering intervals. He also proved that such a phenomenon cannot appear if $f$ is assumed to be $C^2$: every $C^2$-diffeomorphism with irrational rotation number is topologically conjugated to the corresponding irrational rotation.  Note that the conjugating homeomorphisms (or the semiconjugacy) is unique up to composition by a rotation.

However, even for $C^2$ or $C^\infty$ or even analytic diffeomorphorphism with irrational rotation number, the conjugating homeomorphism is \emph{in general} not differentiable. The expression  \emph{in general} here leads to important and deep works, in particular by V. Arnold \cite{Ar}, M. Herman \cite{He} and J.C. Yoccoz \cite{Yo}.  Indeed, for rotation numbers satisfying a diophantine condition, every smooth diffeomorphism is smoothly conjugated to a rotation.
   Later, different  proofs and some generalizations were given  by K. Khanin, Y. Sinai \cite{KhaSi1}, \cite{KhaSi2} and Y. Katznelson, D. Ornstein \cite{KaOr1},\cite{KaOr2}.

We consider here $C^1$-diffeomorphisms. In this class of regularity, no arithmetic condition may ensure a regularity on the conjugacy homeomorphism. Even if we did not find references for this precise statement, it is not doubtful that every irrational rotation number corresponds to infinitely many $C^1$-conjugacy classes. Let us illustrate this different conjugacy classes by distinct behaviors:
\begin{itemize}
\item the $C^1$-centralizer of a diffeomorphism $f$ is the group of diffeomorphisms commuting with $f$. Any diffeomorphism $g$ $C^1$-conjugated to $f$ has a $C^1$-centralizer conjugated to the one of $f$ (by the same diffeomorphism). Therefore, the isomorphism class of the centralizer is a $C^1$-invariant for a $C^1$-conjugacy class: in particular, if $f$ is $C^1$-conjugated to a rotation then its $C^1$-centralizer is isomorphic to $S^1$. There are examples of diffeomorphisms for which the centralizer is  trivial, or some dense subgroup of $\RR$, or much larger than $\RR$ if $f$ admits wandering intervals.
\item  the  asymptotic behavior of the iterates $f^n$ leads also to be an invariant of a $C^1$ conjugacy class: if a $C^1$-diffeomorphism is $C^1$ conjugated to a rotation its derivatives $df^n$ are uniformly bounded for $n\in \ZZ$. However \cite[Theorem B]{BoCrWi}  implies that  any rotation number contains a $C^1$-diffeomorphism for which the sequences $\sup\{df^n(x),df^{-n}(x)\}$, $n\in\ZZ$ is unbounded in any orbit.
\end{itemize}
 All these kind of  properties are invariant under $C^1$- conjugacy, and they show the great variety of $C^1$-behavior of $C^1$-diffeomorphisms with the same irrational rotation number.

 In this paper we consider the space of diffeomorphism having a given irrational rotation number $\alpha\in(\RR\setminus \QQ)/\ZZ$. In his thesis Herman denotes by $F^r_\alpha\subset Diff^r(S^1)$ the closed subset of $C^r$-diffeomorphisms whose rotation number is $\alpha$. He proved several results on $F^r_\alpha$: it is connected and $F^s_\alpha$, for $s>r$ is dense in $F^r_\alpha$ for the $C^r$-topology.  As we said in our comments above, $F^1_\alpha$ always contains many different $C^1$-behavior.  The aim of this paper is to show that these behaviors are indeed equidistributed in $F^1_\alpha$, giving some homogeneity of this space. More precisely:

Given any diffeomorphism $f\in Diff^1(S^1)$ and $r\in\NN$ we denote by $\cC^r(f)$,  its $C^r$-conjugacy class $\{hfh^{-1}, h\in Diff^r(S^1)\}$; notice that  elements in the class $\cC^1(f)$ share the same $C^1$-properties as $f$ (same $C^1$-centralizer, same distorsion properties, etc.). We prove:

\begin{theo}\label{t.conjugacy} Given any $\alpha\in (\RR\setminus \QQ)/\ZZ$ and any $f\in F^1_\alpha$, the $C^1$-conjugacy class $\cC^1(f)$ of $f$ is dense in $F^1_\alpha$ for the $C^1$-topology.
\end{theo}

Approaching the conjugation diffeomorphism $h$ by some smooth diffeomorphism, one gets that $\cC^r(f)$ is also $C^1$-dense, for every $r\in\NN$.

The same kind of question for the rational rotation number case can also be considered. That case  is very related to the question of  conjugacy classes for diffeomorphisms on the segment $[0,1]$: this problem is  solved   in \cite{Fa} which gives complete (and different) answers to  two natural questions:
\begin{itemize}
\item Under what conditions the $C^1$-conjugacy class of a diffeomorphism $f$ of $[0,1]$ contains $g$ in its closure?
\item Under what conditions does it exist a path $h_t$, $t\in[0,1)$ of diffeomorphisms so that $ h_0= Id$ and  $h_tfh_t^{-1}$ tends to $g$ for $t\to 1$?
\end{itemize}

This approach suggest a natural question in our setting:

\begin{ques} Given $f,g\in F^1_\alpha$ does it exist a path $h_t$, $t\in[0,1)$
of diffeomorphisms of $S^1$ so that $ h_0= Id$ and  $h_tfh_t^{-1}$ tends to $g$ for $t\to 1$?
\end{ques}

After announcing our results and the ones of \cite{Fa}, A.  Navas \cite{Na1}
 found a very simple, elementary and clever argument that answers partially to this question;
He showed that

\begin{theo}[Navas] Given any irrational $\alpha$ and  $f\in F^1_\alpha$,  there is   a path $h_t$, $t\in[0,1)$
of diffeomorphisms of $S^1$ so that $ h_0= Id$ and $h_tfh_t^{-1}$ tends to $R_\alpha$ for $t\to 1$.
\end{theo}
Navas argument consists in building the derivatives of the conjugacy $h_t$ as an approximative solution of a cohomological equation, the rotation $R_\alpha$ being characterized in $F^1_\alpha$ by its vanishing logarithm of its derivative. This argument seem not being adaptable for going from $f$ to $g$, when $g$ is not smoothly conjugated to the rotation.

 Notice that a similar result had been proved by Herman in \cite{He} for $C^2$-diffeomorphisms: he proved in that setting that $f$ can be conjugated arbitrarily close to the rotation in the $C^{1 +}$ bounded variations topology.

Given two diffeomorphism $(f_0, g_0) \in Diff^1(S^1) \times Diff^1(S^1)$ and $r\in\NN$ we denote by $\cC^r(f_0,g_0)$,  its $C^r$-conjugacy class $\{(f,g) , f=hf_0h^{-1}, g=hg_0h^{-1} h\in Diff^r(S^1)\}$.
One of our motivation for this paper is the same question for commuting diffeomorphisms
\begin{ques}\label{q.commuting} Given two irrational rotation numbers $\alpha, \beta$, we consider the space of $C^1$ commuting diffeomorphisms $f, g$ with respective rotation numbers $\alpha$ and $\beta$, endowed with the $C^1$ topology.

Are all the $C^1$- conjugacy classes dense  in this space?
\end{ques}

This problem is very related with a famous old question posed by Rosenberg: does it exists pair $(f,g)$ such that the $\ZZ^2$ action is $C^r$-structurally stable?
A positive answer to Question~\ref{q.commuting} would answer negatively to Rosenberg question for $r=1$.
In that direction, Navas \cite{Na1} proved  recently that every conjugacy class contains the pair of rotations $(R_\alpha,R_\beta)$ in its closure.

Notice that, in higher differentiability,  \cite{KlNa} \cite{DeKlNa} provide a generalization of Denjoy theorem for $\ZZ^n$ actions on the circle by $C^{1+\theta}$-diffeomorphisms, where the number $\theta\in(0,1)$ depends on $n$.
For smooth action  J. Moser \cite{Mo}  posed the problem of smooth linearization of commuting circle diffeomorphisms.   In this direction Fayad and Khanin \cite{FaKha} proved that  a finite number of commuting $C^{\infty}$ diffeomorphisms, with simultaneously Diophantine rotation numbers are smoothly conjugated to rotations.

\subsection{Idea of the proof and organization of the paper}

The idea of the proof is very simple. Given $f$ and $g$ with the same irrational rotation number, we want to build a conjugate $hfh^{-1}$ of $f$ arbitrarily $C^1$-close to $g$.
For that, we consider  long orbit segments $x\dots, f^n(x)$ and $y\dots,g^n(y)$ of the same length. They are ordered in the same way on the circle.

 Therefore one may consider a homeomorphism $H$ of the circle so that $H(f^i(x))=g^i(y)$ for $0\leq i\leq n$ and we can choose $H$ being affine on each connected component of the complement  of the orbit segment. If $n$ is large enough and if $f$ and $g$ have dense orbits, the connected component of the complement of these orbit segments are arbitrarily small so that $f$ and $g$ are almost affine on each component, and the derivative on each component is almost the ratio between the component and its image. The same happens for the piecewise $C^1$ homeomorphism $HfH^{-1}$.  Furthermore, up to the components starting at $y$ or at $g^n(y)$ (i.e. the extremities of the orbit segment) the image of a component by $g$ and by $HfH^{-1}$ are the same: as a direct consequence their derivatives are almost equal.
  For the derivatives of $HfH^{-1}$ and $g$ being everywhere (that is, even on the components adjacent to $y$ and $g^n(y)$) almost equal, we show that it is sufficient that the ratio between the length of  adjacent components  to the extremal points $x$, $f^n(x)$ and $y$, $g^n(y)$ are the same for $f$ and for $g$ . This ratio of the lengths of the  components adjacent  to the initial and end point of the orbit segment are called \emph{the initial and final ratio} of $f$ and $g$.

Then, the announced diffeomorphism $h$ is a smoothing of $H$. This is not so easy because the derivative of $H$ can be very different at the right and the left of a singular point, but Proposition~\ref{p.smoothing} solves this difficulty.

Another difficulty comes from the fact that $f$ or $g$ may not have dense orbits, when we deal with $C^1$-diffeomorphisms. The argument can be adapted in that case, once one notices that one may perform a $C^1$ conjugacy so that the distorsion on the wandering interval is arbitrarily small (see Proposition~\ref{p.denjoy}): thus the diffeomorphism is still almost affine on the complement of long orbit segments.

For concluding the proof it remains to show that on can perform small perturbation of $g$ so that its initial and final ratio will be equal to the ones of $f$.  For performing such a perturbation, we would like that the adjacent components  to the extremal points are disjoint from their iterates during a long time, allowing to change their ratio slowly.  This is not always the case. For that, we need to choose carefully the length $n$ of the orbit segments.  We build a sequence of times $k_i$ called \emph{characteristic times} satisfying the property of having a long wandering time.
Lemma~\ref{l.subexp} give a bound of the ratio.  This bound allow us to show at Proposition~\ref{p.caracteristic} that a small perturbation of $g$ at the characteristic times enables us  to get every possible initial and final ratio of $f$, ending the proof.




%% file: Geometrie.tex
\section{Geometry of orbit segments}\label{s.geometry}

In this section we define the fundamental tools of the proof: for every diffeomorphism $f$ with irrational rotation number $\alpha$ we consider orbit segments $x,\dots f^n(x)$, forbidding some exceptional relative position of the first and end point; we call them \emph{adapted segment}.  For this adapted segments we define the initial and final ratio which are the ratio of the lengths of the  components adjacent to $x$ and to $f^n(x)$.

We consider diffeomorphisms $f$ and $g$  with the same irrational rotation number and admitting adapted segments $\{f^i(x)\},\{g^i(y)\}$, $0\leq i\leq n$ of the same lengths $n$ and  with the same initial and final ratio. Then we consider the piecewise affine homeomorphism $H$ sending $f^i(x)$ to $g^i(y)$, $0\leq i\leq n$, which is affine on each component of $S^1\setminus \{f^i(x)\}$.  We show that, if the distorsion of $f$ and $g$ on the components of the complement of the orbits segments is small then the conjugate of $f$ by $H$ is a piecewise $C^1$ homeomorphism whose derivative at each point is close to the one of $g$.

 We notice that this is the case when $f$ and $g$ have dense orbit.  If $f$ has wandering interval, we solve the difficulty is Section~\ref{s.denjoy} by conjugating $f$ to a diffeomorphism with small distorsion on wandering intervals.

 Then we show that we can get a smooth conjugacy of $f$ close to $g$ by smoothing the homeomorphisms $H$ in section~\ref{s.smoothing}.

\subsection{Adapted segments, initial and final ratio, and conjugacy}

Two sequences $x_1,\dots,x_n$ and $y_1,\dots, y_n$ of points of $S^1$ are \emph{similarly ordered} on the circle if there is a orientation preserving homeomorphism $\varphi\colon S^1\to S^1$ with $\varphi(x_i)=y_i$ for  $i\in\{1,\dots,n\}$.

An \emph{orbit segment of length $n$},   of a diffeomorphism $f$ is a sequence $\{x,f(x),\dots,f^n(x)\}$.
Given an orbit segment $\{x,f(x),\dots,f^n(x)\}$; $n\geq 2$, we call \emph{initial and final basic intervals} of  the orbit segment the interval $[a,b]$ and $[c,d]$, respectively,  such that:
\begin{itemize}
\item $a,b,c,d\in \{x,f(x),\dots,f^n(x)\}$
\item $x\in(a,b)$ and $\{x\}=(a,b)\cap \{x,f(x),\dots,f^n(x)\}$
\item $f^n(x)\in(c,d)$ and $\{f^n(x)\}= (c,d)\cap\{x,f(x),\dots,f^n(x)\}$
\end{itemize}

So there is $i,j\in\{1,\dots,n-1\}$ such that $c=f^i(x)$ and $d=f^j(x)$. Notice that $a=f^{n-j}(x)$ and $b=f^{n-i}(x)$.

An orbit segment $\{x,\dots,f^n(x)\}$ is \emph{adapted} if $i\neq 0$, $j\neq 0$, $i+j\neq n-1$

\begin{lemm} If $\{x,\dots,f^n(x)\}$ is  an adapted orbit segment with initial and final basic intervals $[a,b]$ and  $[c,d]$ then the open intervals $(f(c),f(d))$ and $(f^{-1}(a),f^{-1}(b))$ are disjoint from the orbit segment $\{x,\dots,f^n(x)\}$.
\end{lemm}
\begin{demo}  As $i,j$ are different from $n$ by definition, $f(c)$ and $f(d)$ are points of the orbit segments $\{x,\dots,f^n(x)\}$.  Assume (arguing by absurd) $(f(c),f(d))\cap \{x,\dots,f^n(x)\}\neq \emptyset$. That is, there is $k\in\{0,\dots n\}$ with $f^k(x)\in(f(c),f(d))$; then
$f^{k-1}(x)\in(c,d)$.  If $k\neq 0$ this contradics the definition of $(c,d)$. Therefore $(f(c),f(d))\cap \{x,\dots,f^n(x)\}= x$. This means that $i+1=n-j$, that is $i+j=n-1$ contradicting the definition of adapted segment.

This contradiction proves that $(f(c),f(d))\cap \{x,\dots,f^n(x)\}= \emptyset$. The proof of $(f^{-1}(a),f^{-1}(b))\cap \{x,\dots,f^n(x)\}= \emptyset$ is analogous.
\end{demo}

\begin{rema} Given a irrational rotation number $\alpha$, the fact that an orbit segment $\{x,\dots,f^n(x)\}$ is  an adapted orbit segment only depends on the length  $n\geq 0$: more precisely, if $f,g\in F^0_\alpha$ and if $\{x,\dots,f^n(x)\}$ is an adapted orbit segment for $f$, then for every $y\in S^1$, $\{y,\dots,g^n(y)\}$ is an adapted orbit segment.
\end{rema}

Given an adapted segment interval $\{x,\dots,f^n(x)\}$ we call \emph{initial and final ratio} the quotient
$$R_0= \frac{ \ell([a,x])}{\ell([x,b])}\mbox{ and } R_n=\frac{ \ell([c,f^n(x)])}{\ell([f^n(x),d])}.$$

\begin{theo}\label{t.geometry} Let $f$ and $g$ be two diffeomorphisms with irrational rotation number $\alpha$, both $f$ and $g$ with dense orbits.  Assume that, for any
$\varepsilon>0$ there is $\tilde g$ $\varepsilon$-$C^1$-close to $g$, $n\in\NN$ and adapted orbit segments $\{x,\dots, f^n(x)\}$ and $\{y,\cdots,\tilde g^n(y)\}$ similarly ordered and having the same initial and final ratios.

Then there are diffeomorphisms $h_n$ such that $h_n f h_n^{-1}$ tends to $g$ in the $C^1$-topology.
\end{theo}

The diffeomorphism $h_n$ will be obtained as a smoothing of the piecewise affine homeomorphism $H_n$ defined by $H_n(f^k(x))=\tilde g^k(y)$, $k\in\{0,\dots n\}$ and $H_n$ affine in the connected components of the complement of the orbit segment.

\begin{defi} The distorsion $\De(g,I)$ of a diffeomorphism $g$ on some compact interval $I$ is the  maximum of the logarithm of the  quotient of the derivative of $g$ at two point of $I$:
$$\De(g,I)=\max_{x,y\in I} \log (\frac{dg(x)}{dg(y)}).$$
\end{defi}

\begin{lemm}\label{l.affine} Consider diffeomorphisms $f$, $g$ and  a positive number $\varepsilon>0$.  Assume that $f$ and $g$ admits adapted orbit segment
$x\dots,f^n(x)$ and $y,\dots, g^n(y)$ similarly ordered and with the same initial and final ratios.  Assume furthermore that  the distorsion of $f$ and $g$ on each connected component of the complement of the orbit segment is bounded by $\varepsilon$.

Consider the piecewise affine homeomorphism $H$ defined as $H(f^i(x))=g^i(y)$ and $H$ is affine on the connected components of the complement of the orbit segment.

Then:
\begin{itemize}
\item $HfH^{-1}$ is a piecewise $C^1$-diffeomorphism whose derivatives (at the right and at the left) at every point  is close to the derivative of $g$ at the same point. More precisely, $\exp(-2\varepsilon)\leq \frac{d(HfH^{-1})}{d(g)}\leq \exp(2\varepsilon).$

\item $H$ is differentiable at $x$ and $f^n(x)$
\item  $HfH^{-1}$ is $C^1$ up to the segment $\{y,\dots,g^n(y)\}$
\end{itemize}
\end{lemm}
\begin{demo}

 $HfH^{-1}$ is a piecewise $C^1$-diffeomorphism as it is a composition of a piecewise $C^1$-diffeomorphisms.

Notice that $H$ is affine from $[a(f),x]$ to $[a(g),y]$ and from $[x,b(f)]$ to $[y,b(g)]$ and the ratios $\frac{ \ell([a(f),x])}{\ell([x,b(f)])}$ and $\frac{ \ell([a(g),y])}{\ell([y,b(g)])}$ are equal, then $\frac{ \ell([a(f),x])}{\ell([a(g),y])}=\frac{ \ell([x,b(f)])}{\ell([y,b(g)])}$.  This implies that $H$ has the same derivatives at the right and the left sides of $x$, hence is affine in $[a,b]$ (and so smooth at $x$). The proof that $H$ is affine on $[c,d]$ (and so differentiable at $f^n(x)$) is analogous using the final ratios of $f$ and $g$.

Thus $H^{-1}$ is differentiable out of $\{g(y),\dots g^{n-1}(y)\}$.

Furthermore a point $z$ is not singular for $HfH^{-1}$ if $z\notin \{g(y),\dots,g^{n-1}(y)\}$ and $f(H^{-1}(z))\notin \{f(x),\dots, f^{n-1}(x)\}$ that is if $z\notin  \{y,\dots,g^{n-1}(y)\}$.

It remains to compare the derivative of $HfH^{-1}$ with the derivative of $g$. For that, notice that on each connected component $C$ of the complement of $\{y,\dots,g^n(y)\}$ the map is the composition of affine maps with the restriction of $f$ to a connected component of the complement of $\{x,\dots, f^n(x)\}$. Composing with affine maps does not modify the distorsion.  Therefore the distorsion of $HfH^{-1}$ on $C$ is bounded by $ \varepsilon$.
The distorsion of $g$ on $C$ is also bounded by $\varepsilon$.  Furthermore, by construction, $HfH^{-1}(C)=g(C)$. This implies that there is at least a point in $C$ where the derivative of $HfH^{-1}$ and $g$ coincide.  As a consequence
$$\exp (-2\varepsilon)\leq \frac{d(HfH^{-1})}{d(g)}\leq \exp(2\varepsilon).$$

\end{demo}

\begin{rema}\label{r.affine} Notice that the derivatives (at the right and at the left) of $HfH^{-1}$ at every point  is $\varepsilon_0$-close to the derivative of $g$, where $\varepsilon_0=(\exp (2\varepsilon)-1) M< 3 \varepsilon M$ and $M= \sup_ {x \in S^1} |dg(x)|$.

\end{rema}

\begin{prop}\label{p.smoothing} Let $f$ be a $C^1$ diffeomorphism of the circle, $\varepsilon>0$  and $\{x,\dots, f^n(x)\}$ be an adapted orbit segment. Let $H$ be a piecewise affine homeomorphism, smooth out of $\{f(x),\dots, f^{n-1}(x)\}$ such that the right and left derivatives of  $HfH^{-1}$  are $\varepsilon$ close at each point $y\in S^1$.

Then there is a smooth diffeomorphism $h$ arbitrarily $C^0$-close to $H$ and such that the derivative of $hfh^{-1}$ is $2\varepsilon$-close to the right and left derivative of $HfH^{-1}$ at every point.
\end{prop}

We postpone the proof of Proposition~\ref{p.smoothing} to the next section.

\begin{demo}[Proof of Theorem~\ref{t.geometry}] Consider  $f$, $g$, and $\varepsilon$. As the orbits of $f$ and $g$ as assumed to be dense,  we can choose $n$ such that $x,\dots f^n(x)$ and $y,\dots, g^n(y)$ are adapted segment such that the distorsion of $f$ and $g$ on the connected components of the complements of the respective orbit segments is bounded by $\varepsilon/100$.  Choosing $n$ large enough we may approach  $g$ by $\tilde g$  such that  the corresponding orbit segment of $\tilde g$ is still ordered similarly, has the same initial an final ration as the orbit segment of $f$ and the distorsion of $\tilde g$ on the complement of its orbit segment is bounded by $\varepsilon$.

Then Lemma~\ref{l.affine} build a piecewise $C^1$ conjugate of $f$ which is $\varepsilon_0$-$C^1$-close to $\tilde g$ and Proposition~\ref{p.smoothing} allows us  smoothing this piecewise conjugation keeping the $C^1$-proximity to $\tilde g$, hence to $g$.

\end{demo}

\subsection{Smoothing a piecewise linear conjugacy: proof of Proposition~\ref{p.smoothing}}\label{s.smoothing}

We start by linearizing the diffeomorphism $f$ in a neighborhood of the adapted orbit segment:
\begin{lemm} \label{l.linearising} Let $f$ be a diffeomorphism of $S^1$ and $x,\dots,f^n(x)$ be an orbit segment of a non periodic point.
There is a family  of diffeomorphisms $\varphi_t\colon S^1\to S^1, t\in (0,t_0]$ with the following properties
\begin{itemize}
\item $\varphi_t(f^i(x))=f^i(x)$ for every $t$ and $i\in\{0,\dots,n+1\}$,
\item $\varphi_t$ tends to $id|_{S^1}$ for the $C^1$-topology, for $t\to 0$.
\item the derivative of $\varphi_t$ at $f^i(x)$, $i\in\{0,\dots,n+1\}$ is equal to $1$,
\item The restriction of $\varphi_t^{-1} f\varphi_t$ to each segment $[f^i(x)-t,f^i(x)+t]$ is the affine map onto $[f^{i+1}(x)- t\cdot df(f^i(x)), f^{i+1}(x)+ t \cdot df(f^i(x))]$.
\end{itemize}
\end{lemm}
(The proof of this lemma is easy. It is just a change of coordinates, once one notes that all the $f^i(x)$ are distinct.)

For every positive $\alpha,\beta$, we denote by $h_{\alpha,\beta}$ the map define by
\begin{itemize}
\item $x\mapsto \alpha x$ for $x<-1$
\item $x\mapsto \frac{\beta-\alpha}4 x^2+\frac{\beta+\alpha}2 x+\frac{\beta-\alpha}4$ for $x\in[-1,1]$,
\item $x\mapsto \beta x$ for $x>1$
\end{itemize}

An elementary calculation shows that

\begin{lemm} $h_{\alpha,\beta}$ is a $C^1$ diffeomorphism of $\RR$ whose derivative is everywhere contained in $[\alpha,\beta]$.
\end{lemm}

\begin{lemm} \label{l.smoothing} Given positive $\alpha,\beta,\gamma,\delta$ and a point $x\in \RR$ one has
$$\min\{\frac{\alpha}{\gamma},\frac{\beta}{\delta}\} \leq \frac{dh_{\alpha,\beta}(x)}{dh_{\gamma,\delta}(x)}\leq \max\{\frac{\alpha}{\gamma},\frac{\beta}{\delta}\}$$
\end{lemm}
\begin{demo}The proof is straightforward if $x\notin [-1,1]$, since  the maps $h_{\alpha,\beta}$ and $h_{\gamma,\delta}$  are linear with slope $\alpha$ and $\gamma$ (if $x<-1$) or $\beta$ and $\delta$ (if $x>1$).

For $x\in [-1,1]$ one has: $dh_{\alpha,\beta}(x)=\frac{\beta-\alpha}2 x+\frac{\beta+\alpha}2= \frac{1-x}2 \alpha +\frac {1+x}2 \beta$ and $h_{\gamma,\delta}=\frac{1-x}2 \gamma +\frac {1+x}2 \delta$, so that
$$\frac{dh_{\alpha,\beta}(x)}{dh_{\gamma,\delta}(x)}=\frac{\frac{1-x}2 \alpha +\frac {1+x}2 \beta}{\frac{1-x}2 \gamma +\frac {1+x}2 \delta}$$

 The announced inequality now follows immediately from the following (classical) claim:

 \begin{clai}Let $a,b,c,d$ be positive numbers.  Then
 $$\inf\left\{\frac ac,\frac bd\right\}\leq \frac{a+b}{c+d}\leq \max\left\{\frac ac,\frac bd\right\}.$$
 \end{clai}
\begin{demo}[Proof of the claim] Assume $\frac ac\leq \frac bd$ (the converse case is similar). then $a\leq \frac{cb}d$.  Therefore $\frac{a+b}{c+d}\leq \frac{\frac{cb}d+b}{c+d}= \frac{\frac{(c+d)b}d}{c+d}=\frac bd$

This inequality  applied now to $\frac{c+d}{a+b}$ gives now $\frac{c+d}{a+b}\leq\frac ca$ that is
$$\frac ac \leq\frac{a+b}{c+d}\leq \frac bd$$
 which is the announced inequality in that case.
\end{demo}

\end{demo}

If $H$ is a piecewise affine diffeomorphism with a singular point at $x$ and $\alpha,\beta$ be the derivative at the right and the left of $x$, and $\eta>0$ be small enough so that $H$ is affine on $[x-\eta,x]$ and on $[x,x+\eta]$, we denote by $h_{\alpha,\beta,x,\eta}$ the diffeomorphism define on $[x-\eta,x+\eta]$ to $H([x-\eta,x+\eta])$ which is  $B^{-1}\circ h_{\alpha,\beta}\circ A$ where $A$ is the orientation preserving affine diffeomorphism that send $[x-\eta,x+\eta]$ on $[-1,1]$ and $B$ is the orientation preserving affine diffeomorphism sending $H([x-\eta,x+\eta])=[H(x)-\alpha\eta,H(x)+\beta\eta]$ on $[-\alpha, \beta]$.  Notice that

\begin{enumerate}
\item the linear part of $A$ and $B$ coincide, therefore the derivative  $d h_{\alpha,\beta,x,\eta}(z)$ is $d h_{\alpha,\beta}(A(z))$.
\item the derivative of $h_{\alpha,\beta,x,\eta}$ and of $H$ coincide on $x-\eta$ and $x+\eta$.
\item  if $H$ is smooth at $x$, that is $\alpha=\beta$ then $h_{\alpha,\beta,x,\eta}$ coincides with $H$.
\end{enumerate}

\vskip 2mm

We are now ready for proving Proposition~\ref{p.smoothing}:

\vskip 2mm
\begin{demo}[Proof of Proposition~\ref{p.smoothing}]
Let $f$ be a diffeomorphism and $\{x,\dots, f^n(x)\}$ an adapted orbit segment. Let $H$ be a piecewise affine homeomorphism, smooth out of $\{f(x),\dots, f^{n-1}(x)\}$ such that the right and left derivatives of  $HfH^{-1}$  are $\varepsilon$ close at each point $y\in S^1$.

Up to replacing $f$ by a conjugate $\varphi_t^{-1}f\varphi_t$ given by lemma~\ref{l.linearising}, one may assume that there is $t>0$  such that  $f$ is affine in restriction to each interval $[f^i(x)-t, f^i(x)+t]$ for $i\in\{0,\dots,n\}$.

Notice that, for any $\eta>0$ small enough,  and $i\in\{0,\dots, n+1\}$, the interval $f^i([x-\eta,x+\eta])$ is contained in $[f^i(x)-t,f^i(x)+t]$  where $f$ is affine and  $f^i([x-\eta,x+\eta])= [f^i(x)- df^i(x)\cdot \eta, f^i(x)+ df^i(x)\cdot \eta]$.

Let us denote for simplicity:
\begin{itemize}
\item $\eta_i=df^i(x)\cdot \eta$.
\item  $\alpha_i, \beta_i$ are the left and right derivative of $H$ at $f^i(x)$.
\item $A_i\colon [f^i(x)-\eta_i,f^i(x)+\eta_i]\to [-1,1]$ and $B_i\colon [H(f^ix)-\alpha_i\eta_i,H(f^i(x))+\beta_i\eta_i] \to [-\alpha_i, \beta_i]$ are the orientation preserving affine maps.
\end{itemize}

We denote by $h_\eta$ the diffeomorphism of $S^1$ defined as follows:

\begin{itemize}
\item  $h_\eta$ coincide with $H$ out of $\bigcup_{i=1}^{n-1} [f^i(x)-\eta_i, f^i(x)+\eta_i]$
\item $h_\eta= h_{\alpha_i,\beta_i,f^i(x),\eta_i}$ on $[f^i(x)-\eta_i, f^i(x)+\eta_i]$

\end{itemize}

Consider $h_\eta f h_\eta^{-1}$.  For $x \notin H([f^i(x)-\eta_i, f^i(x)+\eta_i])$, $i\in\{0,\dots,n-1\}$ , $h_\eta f h_\eta^{-1}(x) =H f H^{-1}$ so that there is nothing to prove.

Consider $y\in H([f^i(x)-\eta_i, f^i(x)+\eta_i])$.  Then

$$h_\eta f h_\eta^{-1}(y)=h_{\alpha_{i+1},\beta_{i+1},f^{i+1}(x),\eta_{i+1}}\circ f\circ h_{\alpha_i,\beta_i,f^i(x),\eta_i}^{-1}(y)$$

Thus, setting $z=h_{\alpha_i,\beta_i,f^i(x),\eta_i}^{-1}(y)$,  the derivative   is
$$\begin{array}{ccc}
d h_\eta f h_\eta^{-1}(y)&=&df(z)\cdot \frac{dh_{\alpha_{i+1},\beta_{i+1},f^{i+1}(x),\eta_{i+1}}(f(z))}{dh_{\alpha_i,\beta_i,f^i(x),\eta_i}(z)}\\
 &=&  df(z)\cdot \frac{dh_{\alpha_{i+1},\beta_{i+1}}(A_{i+1}(f(z)))}{dh_{\alpha_i,\beta_i}(A_i(z))}
\end{array}
$$

From the fact that $f$ is affine and from the definition of $A_i$ and $A_{i+1}$ one easily check that $A_{i+1} f= A_i$. This implies
$$d h_\eta f h_\eta^{-1}(y)= df(z)\cdot \frac{dh_{\alpha_{i+1},\beta_{i+1}}(A_i(z))}{h_{\alpha_i,\beta_i,}(A_i(z))}.$$

Since $z \in [f^i(x)-\eta_i,f^i(x)+\eta_i]$, one has that

$$d h_\eta f h_\eta^{-1}(y)= df(f^i(x))\cdot \frac{dh_{\alpha_{i+1},\beta_{i+1}}(A_i(z))}{h_{\alpha_i,\beta_i,}(A_i(z))}.$$

According to Lemma~\ref{l.smoothing} one deduces

$$df(f^i(x))\inf\left\{\frac{\alpha_{i+1}}{\alpha_i}, \frac{\beta_{i+1}}{\beta_i} \right\}\leq d h_\eta f  h_\eta^{-1}(y) \leq {df(f^i(x))}\max\left\{\frac{\alpha_{i+1}}{\alpha_i}, \frac{\beta_{i+1}}{\beta_i}\right\}$$

Recall that the derivative of $HfH^{-1}$ is $\frac{\alpha_{i+1}}{\alpha_i}df(f^i(x))$ on $H([f^i(x)-\eta_i,f^i(x)]$  and is $\frac{\beta_{i+1}}{\beta_i}df(f^i(x))$ on $H([f^i(x), f^i(x)+\eta_i]$.

Therefore, the hypothesis on $H$ is that $\left | \frac{\alpha_{i+1}}{\alpha_i}df(f^i(x))-\frac{\beta_{i+1}}{\beta_i}df(f^i(x))\right | <\varepsilon$.

One deduces that $\left | dHfH^{-1}(y)-d h_\eta f  h_\eta^{-1}(y)\right | <2\varepsilon$, as announced.

\end{demo}

%% file: Denjoy.tex
\subsection{Distorsion in wandering intervals for Denjoy counter examples}\label{s.denjoy}

The aim of this section is to prove the following proposition which allows to generalize Theorem~\ref{t.geometry} removing the assumption of dense orbits.
If $f$ is a $C^1$-diffeomorphism with irrational rotation number, we call \emph{(maximal) wandering interval} the closure of each connected component of the complement of the unique minimal set of $f$.

\begin{prop}\label{p.denjoy} Given any diffeomorphism $f\colon S^1\to S^1$ with  irrational rotation number $\alpha$. Then for any $\varepsilon>0$ there is a diffeomorphism  $h$ such that the distorsion of $g=hfh^{-1}$ on each wandering interval $I$ is bounded by $\varepsilon$ .
\end{prop}
\begin{rema}\label{r.denjoy} \begin{enumerate}
\item Let $f$ be a diffeomorphism with an irrational rotation number and $\eta>0$. Then for any point $x$ belonging to the minimal set (that is, $x$ does not belong to any wandering interval) there is $n_1>0$ such that for every $n>n_1$, the closure $I$ of every connected component of the complement of  orbit segment $x,\dots, f^n(x)$ satisfies one of the following possibilities:
\begin{itemize}
\item either the length of $I$ is smaller than $\eta$
\item or there is one wandering interval $J$ contained in $I$ such that the sum of the length of the two components of $I\setminus J$ is smaller than $\eta$.
\end{itemize}

\item Assume now that the distorsion of $f$ is smaller than $\varepsilon/2$ on each wandering interval. Notice that,  due to the uniform continuity of the derivative of $f$, its distorsion on small intervals, is  very small then there is $n_2>0$ such that for any $x$ in the minimal set, for any $n>n_2$ , the distorsion of $f$ on the each connected component of the complement of the orbit segment $x,\dots, f^{n_2}(x)$ is bounded by $\varepsilon$.
\end{enumerate}
\end{rema}

The proof of Proposition~\ref{p.denjoy} is divided in two main parts. We first perturb the derivative by conjugacy inside the orbits of wandering intervals in order to get small distorsion.  Then we will extend the conjugacy on the circle without changing the distorsion inside the wandering intervals.

\begin{lemm}\label{l.wandering} Let $f$ be a $C^1$-diffeomorphism of $S^1$ with an irrational rotation number, and $\varepsilon>0$. Let $[a,b]$ be a maximal wandering interval.  Then there is a family of diffeomorphisms $h_i\colon f^i[[a,b])\to f^i([a,b])$, $i\in\ZZ$ such that :
\begin{itemize}
\item there is $n_0$ such that $h_i=id|_{f^i([a,b])}$ for $|i|\geq n_0$;
\item the distorsion of $h_{i+1}\circ f\circ h_i^{-1}\colon f^i([a,b])\to f^{i+1}([a,b])$ is bounded by $\varepsilon$
\end{itemize}
\end{lemm}

For proving Lemma~\ref{l.wandering} we will use the following lemma:
\begin{lemm} \label{l.fragmentation}Let $\{f_i\}_{i\in\ZZ}$ be a sequence of diffeomorpisms of $[0,1]$ such that $f_n\to id$ for the $C^1$-topology as $n\to\pm\infty$.
Then, there is $\{g_i\}_{i\in\ZZ}$ , arbitrarily $C^1$ close to the identity map, and $n_0$ such that:
\begin{itemize}
\item for $|i|\geq n_0$ one as $g_i=f_i$
\item $g_{n_0}\circ g_{n_0-1} \circ\cdots \circ g_{-n_0+1}\circ g_{-n_0} =f_{n_0}\circ f_{n_0-1} \circ\cdots \circ f_{-n_0+1}\circ f_{-n_0}$
\end{itemize}
\end{lemm}

\begin{demo} Let $n_1>1$ such that $f_n$ and $f_{n+1} f_n$ are $\varepsilon/2$ close to identity for $|n|\geq n_1$.  We fix $g_i=f_i$ for $i<-n_1$.  Consider $F=f_{n_1}\circ \cdots\circ f_{-n_1}$.  A classical elementary result asserts that any orientation preserving diffeomorphism of $[0,1]$ is the product of finitely many diffeomorphisms arbitrarily close to identity. Therefore we can write
$F= g_{m_1} \circ \cdots\circ g_{-n_1}$  with every  $g_i$ $\varepsilon$-$C^1$ close to identity. Up to add several $g_i$ equal to identity, one may assume without lost of generality that $m_1=n_1+k_1$ with $k_1>0$.

Then we write:
\begin{itemize}
\item $g_{m_1+i}=f_{n_1+2i}f_{n_1+2i-1}$ for $i=1, \dots k_1$.
\item $g_i=f_i$ for $i>m_1+k_1=n_1+2k_1$.
\item $n_0>n_1+2k_1$
\end{itemize}
Thus $g_i$ is $\varepsilon$-close to identity for every $i$ and
$$g_{m_1+k_1}\circ \dots\circ g_{m_1}\circ g_{-n_1}= f_{n_1+2k_1}\circ \dots \circ f_{n_1+1}\circ F= f_{n_1+2k_1}\circ\dots \circ f_{-n_1}$$.

As a direct consequence
$g_{n_0}\circ g_{n_0-1} \circ\cdots \circ g_{-n_0+1}\circ g_{-n_0} =f_{n_0}\circ f_{n_0-1} \circ\cdots \circ f_{-n_0+1}\circ f_{-n_0}$
concluding the proof of the lemma.
\end{demo}

\begin{demo}[Proof of Lemma~\ref{l.wandering}]Let $\varphi_i\colon f^i([a,b])\to [0,1]$ be the orientation preserving affine diffeomophism.  Then we write $f_i=\varphi_{i+1} f\varphi_i^{-1}$.  Notice that for $i$ big, length of $f^i([a,b])$ is small, hence the distorsion of $f$ on $f^i([a,b])$ implies the $C^1$ distance of $f_i$ to identity.

Therefore, the sequence $f_i$ satisfies the hypothesis of Lemma~\ref{l.fragmentation}. Consider $n_0>0$ and  $g_i$  the sequence of diffeomorphisms given by lemma~\ref{l.fragmentation} such that $g_i$ is $\frac{\varepsilon}{4}$  $C^1$ close to identity.  In particular $f_i=g_i$ for $i<-n_0$.

One set
\begin{itemize}

\item  $h_i=id$ for $i<-n_0-1$.
\item $h_i=\varphi_i^{-1}\circ g_{i-1}\circ\cdots\circ g_{-n_0-1}\circ f_{-n_0-1}^{-1} \circ\cdots\circ f_{i-1}^{-1}\circ \varphi_i$
\end{itemize}

By definition of the $g_i$ one can check that $h_i=id$ for $i>n_0$.

Furthermore
$$\begin{array}{ccl}
h_{i+1} f h_i^{-1}&=&\varphi_{i+1}^{-1}\circ \\
&& g_{i}\circ\cdots\circ g_{-n_0-1}\circ f_{-n_0-1}^{-1} \circ\cdots\circ f_{i}^{-1}\circ \\
&& \varphi_{i+1}\circ f\circ\varphi_i^{-1}\circ\\
&&   f_{i-1}\circ\cdots\circ f_{-n_0-1}\circ g_{-n_0-1}^{-1} \circ\cdots\circ g_{i-1}^{-1}\\
&&\circ \varphi_i \\
&=&\varphi_{i+1}^{-1}\circ \\
&& g_{i}\circ\cdots\circ g_{-n_0-1}\circ f_{-n_0-1}^{-1} \circ\cdots\circ f_{i}^{-1}\circ \\
&& f_i\circ\\
&&   f_{i-1}\circ\cdots\circ f_{-n_0-1}\circ g_{-n_0-1}^{-1} \circ\cdots\circ g_{i-1}^{-1}\\
&&\circ \varphi_i \\
&=&\varphi_{i+1}^{-1}\circ g_{i}\circ \varphi_i \\
\end{array}$$

Note that as $g_i$ is $\frac{\varepsilon}{4}$-$C^1$ close to the identity map one has
$|\frac{dg(x)}{dg(y)}-1|<\frac{2\varepsilon}{3}$.

As $\varphi_i$ and $\varphi_{i+1}$ are affine  one gets that the distorsion of $h_{i+1} f h_i^{-1}$ on $f^i([a,b])$ is bounded by $ \log (1+\frac{2\varepsilon}{3})$, therefore it is smaller than $\varepsilon$.
\end{demo}

Next lemma ensures that one can extend the conjugacy, defined inside the wandering interval by lemma~\ref{l.wandering}, on the whole circle without changing the distorsion.

First notice that  due to the uniform continuity of the derivative of $f$, its distorsion is smaller than $\varepsilon/2$ on every small enough interval. Therefore, we can choose  $N, k$ and finitely many orbits segments of (maximal)wandering intervals $[f^{-N}(a_i), f^{-N}(b_i)],\dots [a_i,b_i],\dots ,[f^{-N}(a_i), f^{-N}(b_i)]$, $i\in\{1,\dots k\}$, such that
\begin{itemize}
\item The orbits of the $[a_i,b_i]$ are pairwise distinct
\item for any $n$ with $|n|>N$ the distorsion of $f$ on $f^n([a_i,b_i)$ is smaller than $\frac \varepsilon 2$
\item for any wandering interval $[a,b]$ whose orbit is distinct from the $[a_i,b_i]$ the distorsion is bounded by $\varepsilon/2$ on each $f^i([a,b])$, $i\in\ZZ$.
\end{itemize}

For every $t\in\ZZ$ we chose a diffeomorphism $h_{i,t}\colon [f^t(a_i),f^t(b_i)]\to[f^t(a_i),f^t(b_i)]$ given by lemma~\ref{l.wandering}, so that
\begin{itemize}
\item there is $n_0$ such that $h_{i,t}=id|_{f^t([a_i,b_i])}$ for $|t|\geq n_0$ and any $i\in\{1,\dots, k\}$;
\item the distorsion of $h_{i,t+1}\circ f\circ h_{i,t}^{-1}\colon f^t([a_i,b_i])\to f^{t+1}([a_i,b_i])$ is bounded by $\varepsilon$
\end{itemize}

Notice that the segments $[f^t(a_i),f^t(b_i)]$, $i\in\{1,\dots, k\}$, $t\in\{-n_0,\dots, n_0\}$ are finitely many compact disjoint segments.

\begin{lemm}\label{l.denjoy} There is a diffeomorphism $h$ of $S^1$, such that $h$ coincides with the $h_{i,t}$ for $i\in\{1,\dots, k\}$, $t\in\{-m,\dots, m\}$, and the derivative of $h$ is constant in every wandering interval distinct from the $[f^t(a_i),f^t(b_i)]$, $i\in\{1,\dots, k\}$, $t\in\{-m,\dots, m\}$.
\end{lemm}
\begin{demo} Consider the closure $I=[f^{t_1}(b_i),f^{t_2}(a_j)]$ of a connected component of 

$S^1\setminus \bigcup_{i\in\{1,\dots, k\}, t\in\{-m,\dots, m\}}[f^t(a_i),f^t(b_i)]$.

We build a continuous map $\varphi_I\colon I\to (0,+\infty) $ such that
\begin{itemize}
\item $\varphi_I$ constant on each wandering interval in the interior of $I$
\item $\varphi_I(f^{t_1}(b_i)=d h_{i,t_1}(f^{t_1}(b_i))$ and $\varphi_I(f^{t_2}(a_j)= d h_{j,t_2}(f^{t_2}(a_j))$
\item $\int_I \varphi_I=\ell(I)=|f^{t_1}(b_i)-f^{t_2}(a_j)|$
\end{itemize}
 The existence of such a function is the usual construction of Lebesgue devil stairs.

Then we define $h|_I$ by $h(x)=f^{t_1}(b_i)+ \int_{f^{t_1}(b_i)}^x \varphi_I$.  Its a diffeomorphism of $I$ whose derivative coincides with the one of $h_{i,t_1}$ on $f^{t_1}(b_i)$ and with the one of $h_{j,t_2}$ on $f^{t_2}(a_j)$.

The concatenation for this diffeomorphisms is the announced diffeomorphisms $h$.

\end{demo}


\begin{demo}[Proof of Proposition~\ref{p.denjoy}] Consider the diffeomorphism $h$ given by Lemma~\ref{l.denjoy}.  Then  $hfh^{-1}$ coincides with $h_{i,t+1}fh_{i,t}$ on the $f^t([a_i,b_i])$ for $|t|<n_0$; therefore, the distorsion is bounded by $\varepsilon$.  On the other wandering intervals, $h$ is affine so that the conjugacy does not affect the distorsion, which  was bounded by $\frac{\varepsilon}{2}$ by definition of the $[a_i,b_i]$ and $n_0$.
\end{demo}

%% file: Prova.tex
\section{Proof of the main result}

The aim of this section is to prove  Theorem~\ref{t.conjugacy}
assuming  Theorem~\ref{t.caracteristic} which explains that one can change the initial and final ratio by arbitrarily small perturbations if one choose adapted segment of a specific length.

\subsection{Perturbing the initial and final ratio at \emph{characteristic times}}

Our main technical result is

\begin{theo}\label{t.caracteristic} Given any $\alpha\in\RR\setminus \QQ$ there is a  (strictly increasing) sequence $\{k_i\}\subset \NN$  with the following properties:
\begin{itemize}
\item Given any $C^1$-diffeomorphisms $f$, $g$ with rotation number $\alpha$,
\item given any points $x,y\in S^1$
\item given any $\varepsilon>0$
\end{itemize}

Then,
\begin{itemize}
\item for any $i$ the orbits segments $\{x,\dots, f^{k_i}(x)\}$ and $\{y,\dots, g^{k_i}(y)\}$ are adapted segments,
\item there is $i_0$ such that for every $i\geq i_0$ there is a $C^1$-diffeomorphism $g_i$ such that
\begin{itemize}
\item $g_i$ is $\varepsilon$ $C^1$-close to $g$
\item $\{y,\dots, g_i^{k_i}(y)\}$ is an adapted segment of $g_i$ ordered on $S^1$ as $\{x,\dots, f^{k_i}(x)\}$ and $\{y,\dots, g^{k_i}(y)\}$;
\item The initial and final ratio of $g_i$ on  the adapted segment $\{y,\dots, g_i^{k_i}(y)\}$ are the same as the ones of $f$ on the adapted segment $\{x,\dots, f^{k_i}(x)\}$.
\end{itemize}
\end{itemize}

\end{theo}

In section~\ref{s.characteristic} we will build the sequence $\{k_i\}$ called \emph{characteristic times}, and Section~\ref{s.perturbation} will be dedicated to the proof of Theorem~\ref{t.caracteristic}.

The aim of this section is to show that Theorem~\ref{t.caracteristic} together with Proposition~\ref{p.denjoy},  Lemma~\ref{l.affine} and Proposition~\ref{p.smoothing} imply Theorem~ \ref{t.conjugacy}.

\subsection{Proof of Theorem~\ref{t.conjugacy}}

Let $f, g\in F^1_\alpha$, and $\varepsilon>0$.  We have to prove that there is  a diffeomorphism $h$ of $S^1$ such that $hfh^{-1}$ is $\varepsilon$-close to $g$.

According to \cite[Proposition 4.4.2]{He} the set $F^r_\alpha$ of $C^r$-diffeomorphisms with rotation number $\alpha$ is $C^s$-dense in $F^s_\alpha$ for any $s\leq r$. In particular, $F^2_\alpha$ is $C^1$-dense in $F^1_\alpha$.  Thus there is a $C^2$-diffeomorphism $g_0$ with rotation number $\alpha$ and $\varepsilon/2$-close to $g$. In other word, up to change $\varepsilon$ by $\varepsilon/2$ and $g$ by $g_0$, we may assume without lost of generality that $g$ is $C^2$.

According to Proposition~\ref{p.denjoy} $f$ is $C^1$ conjugated to $f_0=h_0fh_0^{-1}$ such that the distorsion of $f_0$ on each wandering interval is bounded by $\frac\varepsilon{48M}$, where $M$ is an upper bound for $dg$. Therefore, according to Remark~\ref{r.denjoy}, for any sufficiently large orbit segment associated to a point $x$ in the minimal set of $f_0$,
the distorsion of $f_0$ on each connected component of the complement of the orbit segment will be bounded by $\frac\varepsilon{24M}$.

Thus, we chose $x$ in the minimal set of $f$ and $x_0=h_0(x)$ is on the minimal of $f_0$.  We chose a sufficiently large characteristic time $k_i$ so that, according to Theorem~\ref{t.caracteristic},  $g$ admits an $\frac {\varepsilon} {200}$ $C^1$ perturbation $g_1$ for which
\begin{itemize}
\item the orbit segment $\{0,\dots g_1^{k_i}(0)\}$ is adapted and is ordered as $\{0,\dots g^{k_i}(0)\}$
\item the initial and final ratio associated to  the orbit segment $\{0,\dots g_1^{k_i}(0)\}$ are the same as the ones of $f_0$ on $\{x_0,\dots f_0^{k_i}(x_0)\}$.
\item the distorsion of $g_1$ on each connected component of the complement of the orbit segment $\{0,\dots g_1^{k_i}(0)\}$ will be  be bounded by $\frac\varepsilon{24M}$, since
  $g_1$ was chosen $C^1$- close to $g$.
\end{itemize}

Now Lemma~\ref{l.affine} build a piecewise linear conjugacy $H$ so that $Hf_0H^{-1}$ satisfies (see Remark~\ref{r.affine}) that $|dHf_0H^{-1}-dg|<\frac{\varepsilon}{8}$.

Finally Proposition~\ref{p.smoothing} ensures the existence of a diffeomorphism $h$ for which $|dHf_0H^{-1}-dhf_0h^{-1}|<\frac\varepsilon 2$.
One gets that $hh_0 f(hh_0)^{-1}$ is $\varepsilon$ $C^1$-close to $g$, concluding the proof.

\hfill$\Box$\medskip

It remains now to prove Theorem~\ref{t.caracteristic} for concluding the proof.

%% file: Characteristic.tex
\section{Characteristic times}\label{s.caracteristic}
 Section~\ref{s.geometry} shows that the main point for conjugating  a diffeomorphism $f$ in order to be $C^1$ close to $g$ is to control the initial and final ratio of adapted segments.

In this section, we will choose specific adapted segment that we will call \emph{characteristic segments}. This orbit segments will be choosen for the rotation of $R_\alpha$, $\alpha\in\RR\setminus\QQ$, and we will then control the ratio of characteristic segments for the diffeomorphisms $f,g$ in $F^i_\alpha$, $i=1,2$, and for their $C^1$ perturbations.

The idea is very simple:   consider the closest return close to $0$ of  the  $n.\alpha$, and consider the time  $n-1$ just before it; then the ratio  of the corresponding  segments for the rotation $R_\alpha$ is uniformly bounded, between $1/2$ and $2$. Now we will extract a subsequence (called \emph{characteristic time}) for which the union of the two segments adjacent to $0$ will have a large number of disjoint successive iterates.
This long wandering time will allow us to modify as we want this ratio by a $C^1$-perturbation.

As we need to control the complete geometry of the orbit segment until the closest return, we will first reconstruct the sequence of this closest return times, paying attention to the wandering time of union of the segments adjacent to $0$.

\subsection{Ordering the orbit segments of rotations}
\label{s.characteristic}

Let $\alpha\in\RR\setminus \QQ$ be an irrational number and $R_\alpha$ the rotation $x\mapsto x+\alpha$ on $S^1=\RR/\ZZ$.

Every orbit segment $\{x,\dots, R_\alpha^n(x)\}$ of the rotation with length $n$  is the image by the isometry $R_x$ of the corresponding orbit segment stating at $0$.  We consider therefore the orbit segments  $\{0,\alpha,\dots, n\alpha\}$.

We consider the points $-\frac12<-a_n<0<b_n<\frac12$ which are adjacent to $0$ in this orbit segment. We define $r_n,s_n\in\{1,\dots, n\}$ by  $-a_n= r_n\alpha$ and $b_n=s_n\alpha$. Notice that $(n-s_n)\alpha<n\alpha<(n-r_n)\alpha $ are the adjacent points to $n\alpha$.

 Note that as $a_n \rightarrow 0$, and
 $b_n \rightarrow 0$ when  $n \rightarrow \infty$ then $r_n \rightarrow \infty$ and $s_n \rightarrow \infty$.

The following lemma is very elementary and classical.
\begin{lemm}\label{l.classic} \begin{enumerate}
\item  The length of each connected component of the complement of the orbit segment belong to $\{a_n, b_n,a_n+b_n\}$.
\item $$r_n+s_n\neq n+1\Longleftrightarrow R_\alpha(n\alpha)= (n+1) \alpha \notin [-a_n,b_n].$$ In that case:
\begin{itemize}
\item  $(s_n-1)\alpha,(r_n-1)\alpha$ are adjacent and the length of the component $((s_n-1)\alpha,(r_n-1)\alpha)$ is  $a_n+b_n$.
\item $a_{n+1}=a_n$ and $b_{n+1}=b_n$, $r_{n+1}=r_n$, $s_{n+1}=s_n$.
\end{itemize}

\item If $r_n+s_n = n+1$ then the image by $R_\alpha$ of the segment $[(n-s_n)\alpha,(n-r_n)\alpha] $ is the segment $[r_n\alpha, s_n\alpha]=[-a_n,b_n]$.  One deduces that the length of each connected component of the complement of the orbit segment belong to $\{a_n, b_n\}$. Furthermore:
\begin{enumerate}
\item Assume $a_n>b_n$ then $(n+1)\alpha\in(-a_n,0)$ and
$$\left\{\begin{array}{ccl}
a_{n+1}&=&a_n-b_n\\
b_{n+1}&=&b_n\\
r_{n+1}&=&n+1=r_n + s_n\\
s_{n+1}&=&s_n\\
\end{array}\right.
$$
\item Assume $b_n>a_n$  then $(n+1)\alpha\in(0,b_n)$ and
$$\left\{\begin{array}{ccl}
a_{n+1}&=&a_n\\
b_{n+1}&=&b_n-a_n\\
r_{n+1}&=&r_n\\
s_{n+1}&=&n+1= r_n + s_n\\
\end{array}\right.
$$
\end{enumerate}
\end{enumerate}
\end{lemm}




Let $n_i$ denote the sequence of number such that $(n_i+1)\alpha \in[-a_{n_i},b_{n_i}]$.
Notice that:
\begin{itemize}
\item if $a_{n_i}>b_{n_i}$ then $a_{n_{i+1}}=a_{n_i+1}=a_{n_i}-b_{n_i}$ and $b_{n_{i+1}}=b_{n_i+1}=b_{n_i}$.
\item if $a_{n_i}<b_{n_i}$ then $a_{n_{i+1}}=a_{n_i+1}=a_{n_i}$ and $b_{n_{i+1}}=b_{n_i+1}=b_{n_i}-a_{n_i}$.
\end{itemize}

One deduces:

\begin{lemm} There is a subsequence $n_{i_j}$ of $n_i$ such that $\frac{a_{n_{i_j}}}{b_{n_{i_j}}}\in[\frac12,2]$.
\end{lemm}
\begin{demo}Assume $a_{n_i}>b_{n_i}$.  Therefore $a_{n_{i+1}}=a_{n_i}-b_{n_i}$ and $b_{n_{i+1}} = b_{n_i}$.
If $a_{n_{i+1}}<b_{n_{i+1}}$ this means $b_{n_i}<a_{n_i}<2b_{n_i}$ so that $n_i$ belongs to the announced sequence.

Otherwise, $a_{n_{i+1}}>b_{n_{i+1}}$ and $a_{n_{i+2}}=a_{n_{i+1}}-b_{n_{i+1}}<a_{n_i}$;
if  $a_{n_{i+1}}-b_{n_{i+1}}<b_{n_i}=b_{n_i+1}=b_{n_i+2}$ we are done; otherwise  we continue the process so that there is $k$ such that
$a_{n_{i+k}}>b_{n_{i+k}}=b_{n_i}$ but $a_{n_{i+k+1}}<b_{n_{i+k+1}}$: then $n_{i+k}$ belongs to the announced sequence.

The case $a_{n_i}<b_{n_i}$ is analogous.  Thus we have shown that the announced sequence contains numbers greater than any of the $n_i$, allowing to define the $n_{i_j}$ by induction.
\end{demo}

\begin{rema}  The points $(n_{i_j}+1)\alpha$ are the sequence of \emph{closest return} to $0$.
\end{rema}

\subsubsection{Wandering time}
Consider $n>0$, the orbit segment $0,\dots, n\alpha$ and the numbers $r_n\alpha, s_n\alpha$, so that $0$ is the unique point of the segment in the open interval $I_n=(r_n\alpha, s_n\alpha)$ and $n$ is the unique point of the segment in the open interval $J_n=((n-s_n)\alpha, (n-r_n)\alpha)$.

We call \emph{wandering time $w(n)$} the largest integer $w$ such that the intervals $I_n, R_\alpha(I_n),\dots, R_\alpha^w(I_n),$ $ R_\alpha^{-w} J_n,\dots,J_n$ are pairwise disjoint.

\begin{lemm} The wandering time $w(n)$ is
$$w(n)=\inf{[\frac{n-r_n-1}2],[\frac{n-s_n-1}2]}$$
where $[. ]$ denotes the integer part.
\end{lemm}
\begin{demo} The intervals $(r_n\alpha,0),\dots,(n\alpha,(n-r_n)\alpha), (0,s_n\alpha),\dots, ((n-s_n)\alpha,n\alpha)$ are  pairwise disjoint.
We just break this family in two families of equal length.
\end{demo}

Recall that,  the times $n_{i_j}$ (before the closest return) are characterized by $r_{n}+s_{n}= n+1$ and that  the sequences $r_n$ and $s_n$ go to infinity. Furthermore
$$w(n_{i_j})=\inf \{{[\frac{s(n_{i_j})}2]-1,[\frac{r(n_{i_j})}2]-1}\}.$$

Lemma~\ref{l.timeratio} and Corollary~\ref{c.timeratio} below give a lower bound for the wandering time:

\begin{lemm} \label{l.timeratio} Exists a strictly increasing sequence $N_i$ of integer such that $N_i+1$ is the closest return time (that is, $N_i$ is a subsequence of the $n_{i_j}$) and
\begin{itemize}
\item either $a_{N_i}<b_{N_i}$ and   $ r({N_i})\leq 2 s({N_i})$
\item or $a_{N_i}>b_{N_i}$ and   $s({N_i})\leq 2 r({N_i})$.
\end{itemize}
\end{lemm}
\begin{coro}\label{c.timeratio} With the notations of Lemma~\ref{l.timeratio}, for every $i$
\begin{itemize}
\item either $a_{N_i}<b_{N_i}$ and   $ \left[\frac{r(N_i)}4\right] -1 \leq w(N_i)$
\item or $a_{N_i} > b_{N_i}$ and   $ \left[\frac{s(N_i)}4\right] -1 \leq w(N_i)$.
\end{itemize}
\end{coro}

\begin{demo} Denote by $\cN(\alpha)$ the subset of  $\{n_{i_j}\}_{j\in\NN}$ satisfying  either $a_{n_{i_j}}<b_{n_{i_j}}$ and   $ r_{n_{i_j}}\leq 2 s_{n_{i_j}}$ or $a_{n_{i_j}}>b_{n_{i_j}}$ and   $s_{n_{i_j}}\leq 2 r_{n_{i_j}}$.  We have to prove that $\cN(\alpha)$ is infinite, for every $\alpha\in\RR\setminus \QQ$.

Most of the $\alpha$ are solved by the following claim:
\begin{clai} If $i_{j}-i_{j-1}\geq 2$ then
\begin{itemize}
\item either $a(n_{i_j})<b(n_{i_j})$ and $r(n_{i_j})<s(n_{i_j})$
\item or $a(n_{i_j})>b(n_{i_j})$ and $r(n_{i_j})>s(n_{i_j})$
\end{itemize}
so that in both cases  $i_j\in\cN(\alpha)$
\end{clai}
\begin{demo} Assume for instance $a(n_{i_j})<b(n_{i_j})$, the other case is identical.

 By the choice  of the $n_{i_j}$, one has $a(n)<b(n)$ for every $n_{i_{j-1}}<n\leq n_{i_j}$.

Furthermore, by lemma~\ref{l.classic}
$$\begin{array}{lcl}
r(n_{i_{j-1}}+1)&= &r(n_{i_{j-1}})+s(n_{i_{j-1}})\\
s(n_{i_{j-1}}+1)&= &s(n_{i_{j-1}})
\end{array}
$$
Furthermore,  $r(n_{i_{j-1}+1})=r(n_{i_{j-1}}+1)$ and $s(n_{i_{j-1}+1})=s(n_{i_{j-1}}+1)$

Then by lemma~\ref{l.classic}, for every $0\leq k \leq i_j-i_{j-1}$ one has

$$
\begin{array}{lcl}
r(n_{i_{j-1}+k})&= &r(n_{i_{j-1}}+1)\\
s(n_{i_{j-1}+k})&=&s(n_{i_{j-1}}+1) + (k-1)r(n_{i_{j-1}}+1)\\
\end{array}
$$

In particular,
$$
\begin{array}{lcl}
r(n_{i_{j}})&= &r(n_{i_{j-1}}+1)\\
s(n_{i_{j}})&=&s(n_{i_{j-1}}+1) + (i_j-i_{j-1}-1)r(n_{i_{j-1}}+1)\\
\end{array}
$$
so that $r(n_{i_j})<s(n_{i_j})<2s(n_{i_j})$ as announced.
\end{demo}

Consider now $\alpha$ such that there is $j_0$ so that  every $j\geq j_0$ does not satisfy the conclusion of the claim.
This implies that $i_{j_0+k}= i_{j_0}+k$ for every positive $k$.

Assume for instance $a(n_{i_{j_0}})<b(n_{i_{j_0}})$.

Therefore, for every $k> 0$ one has
$$\begin{array}{lcl}
r(n_{i_{j_0}}+2k)&= &r(n_{i_{j_0}}+2k-1)+ s(n_{i_{j_0}}+2k-1)\\
s(n_{i_{j_0}}+2k)&= &s(n_{i_{j_0}}+2k-1)\\
\end{array}
$$

$$\begin{array}{lcl}
r(n_{i_{j_0}}+2k+1)&= &r(n_{i_{j_0}}+2k)\\
s(n_{i_{j_0}}+2k+1)&= &r(n_{i_{j_0}}+2k)+ s(n_{i_{j_0}}+2k)\\
\end{array}
$$

In particular  $r(n_{i_{j_0+2}})=r(n_{i_{j_0+1}})+s(n_{i_{j_0+1}})=2r(n_{i_{j_0}})+s(n_{i_{j_0}})$ and $s(n_{i_{j_0+2}})=s(n_{i_{j_0+1}})=r(n_{i_{j_0}})+s(n_{i_{j_0}})$  so that

$$r(n_{i_{j_0+2}})<2s(n_{i_{j_0+2}})$$
This proves that $n_{i_{j_0+2}}\in \cN(\alpha)$, and ends the proof of the lemma.
\end{demo}

The $N_i$ are almost the announced characteristic times. The unique defect is that the orbit segments $x,\dots, R_\alpha^{N_i}(x)$ are not adapted. We define $k_i=N_i-1$ and we call them the characteristic times. We denote $w_i=w(k_i)$, the wandering time of the characteristic segment.
 Then
\begin{itemize}
\item The orbit segments $x,,\dots, R_\alpha^{k_i}(x)$ are adapted segments,
\item The initial and final ratio of the rotation $R_\alpha$ belong to $[\frac 12,2)$
\item
\begin{itemize}
\item either $a_{k_i}<b_{k_i}$ and   $ \left[\frac{r(k_i)}4\right] -1 \leq w(k_i)$
\item or $a_{k_i} > b_{k_i}$ and   $ \left[\frac{s(k_i)}4\right] -1 \leq w(k_i)$.
\end{itemize}

\end{itemize}





\subsection{Geometry of the characteristic segment for $C^1$ diffeomorphisms $f\in F^1_{\alpha}$}

Let $f$ be a $C^1$-diffeomorphism with an irrational rotation number  $\alpha$. Classical results assert that $f$ is uniquely ergodic, that is adimits a unique invariant measure. The Lyapunov exponent of this measure is zero. This implies:
\begin{lemm} For any $\lambda>1$ there is $n_{\lambda}>0$ such that for any $n> n_\lambda$ and any $x\in S^1$ one has
$$df^n(x)\in[\lambda^{-n},\lambda^n].$$
\end{lemm}

One check easily:
\begin{coro}\label{c.bounds} Let $x,y\in S^1$ such  that there is $n$ with $|n|>n_\lambda$ with  and $x<f^n(x)<f^{2n}(x)<y<f^{3n}(x)$.  Then
$$\frac{|x-f^n(x)|}{\lambda^{|n|}}<|f^n(x)-y|<(\lambda^{|n|}+\lambda^{|2n|})|x-f^n(x)|$$
\end{coro}

\begin{lemm}\label{l.subexp} Let $f\in F^1_\alpha$ be a $C^1$-diffeomorphism with irrational rotation number $\alpha$, and let $k_i$ be the sequence of numbers associated to $\alpha$ defined after Lemma~\ref{l.timeratio}.
Then, given any $\lambda>1$ there is $i(\lambda)$ so that for every $i\geq i(\lambda)$  and any $x\in S^1$ the initial and final ratio of the characteristic segments $x,\dots, f^{k_i}(x)$ belong to $[\lambda^{-w_i},\lambda^{w_i}]$.
\end{lemm}
\begin{demo} Consider the corresponding orbit segment $0,\dots k_i\alpha$ for the rotation $R_\alpha$ and assume, for instance that $a(k_i)<b(k_i)$. By the definition of the characteristic times, we also have $b(k_i)<2a(k_i)$.

This means $R_\alpha^{r(k_i)}(0)<0<R_\alpha^{-r(k_i)}(0)<R_\alpha^{s(k_i)}(0)< R_\alpha^{-2r(k_i)}(0)$

As the order of the orbits just depends of the rotation number, one gets

$f^{r(k_i)}(x)<x<f^{-r(k_i)}(x)<f^{s(k_i)}(x)< f^{-2r(k_i)}(x)$.

Given any $\lambda_1>1$ and $i$ so that $r(k_i)>n_{\lambda_1}$,  one deduce the following bounds  from Corollary~\ref{c.bounds} applied to $f^{r(k_i)}(x)$ and $f^{-r(k_i)}$:

$$\frac{|x-f^{r(k_i)}(x)|}{\lambda_1^{|r(k_i)|}}<|f^{s(k_i)}(x)-x|<(\lambda_1^{|r(k_i)|}+\lambda_1^{|2r(k_i)|})|x-f^{r(k_i)}(x)|$$

Thus the initial ration belongs to $[\lambda_1^{-r(k_i)}, \lambda_1^{|r(k_i)|}+\lambda_1^{|2r(k_i)|}]$

Recall that, for characteristic time $k_i$ for which $a(k_i)<b(k_i)$ one has $ \left[\frac{r(k_i)}4\right] -1 \leq w(k_i)$.  Thus, for concluding the proof of the lemma,  it is enough to choose $\lambda_1$ so that $\lambda_1^{|r(k_i)|}+\lambda_1^{|2r(k_i)|}<\lambda^{w(k_i)}$ for every $k_i>n_{\lambda_1}$.

This give us the announce bounds for the initial ratio, the final ratio is obtained similarly.
\end{demo}

We can restate Lemma~\ref{l.subexp} as follows:
\begin{rema} \label{r.subexp}
There is a sequence $\lambda_i>1$ tending to $1$ as $n\to\infty$, so that for every $i$  and any $x\in S^1$ the initial and final ratio of the characteristic segments $x,\dots, f^{k_i}(x)$ belong to $[\lambda_i^{-w_i},\lambda_i^{w_i}]$.
\end{rema}

%% file: Perturbations.tex
\section{Perturbations}

\label{s.perturbation}

The aim of this section is to prove Theorem~\ref{t.caracteristic} using the characteristic time $\{ k_i\}$ as the announced sequence. Let $w_i$ be the corresponding wandering times.
In the statement of the theorem, the diffeomorphism $f$ appears only by its initial and final ratio. Let us recall that, according to Lemma~\ref{l.subexp} and Remark~\ref{r.subexp} this ratio are bounded in an interval $[\lambda_i^{-w_i},\lambda_i^{w_i}]$ where the sequence $\lambda_i>1$ tends to $1$; these sequence depends on $f$. Finally, the theorem announces an $\varepsilon$ perturbation of $g$.  Therefore, up to shrink $\varepsilon$ if necessary, one may assume that $g$ is a $C^2$ diffeomorphism.

Let us restate the theorem:

\begin{prop}\label{p.caracteristic}
\begin{itemize}
\item Given any $\alpha\in\RR\setminus \QQ$,  $\{k_i\}\subset \NN$  its characteristic times, and $w_i$ the corresponding wandering times.
\item Given any sequence $\lambda_i>1$ tending to $1$,
\item Given any sequences $\rho^-_i,\rho^+_i\in [\lambda_i^{-w_i},\lambda_i^{w_i}]$
\item Given any $C^2$-diffeomorphism  $g$ with rotation number $\alpha$,
\item given any point $y\in S^1$
\item given any $\varepsilon>0$
\end{itemize}

Then, there is $i_0$ such that for every $i\geq i_0$ there is a $C^1$-diffeomorphism $g_i$ such that
\begin{itemize}
\item $g_i$ is $\varepsilon$ $C^1$-close to $g$
\item $\{y,\dots, g_i^{k_i}(y)\}$ is an adapted segment of $g_i$ ordered on $S^1$ as $\{y,\dots, g^{k_i}(y)\}$;
\item The initial and final ration of $g_i$ on  the adapted segment $\{y,\dots, g_i^{k_i}(y)\}$ are $\rho^-_i$ and $\rho^+_i$, respectively.
\end{itemize}

\end{prop}
\subsection{Rescaling the statement of Proposition~\ref{p.caracteristic} on the segment $[0,1]$}

Let $I_i=I_i(g)$ and $J_i=J_i(g)$ denote the intervals $(g^{r(k_i)}(y),g^{s(k_i)}(y))$ and $(g^{k_i-s(k_i)}(y),g^{k_i-r(k_i)}(y))$, containing $y$ and $g^{k_i}(x)$, respectively.
By definition of the wandering times, the intervals \\
$I_i(g)$, $g(I_i),\dots$,$ g^{w_i}(I_i)$, $ g^{-w_i} (J_i),\dots,J_i$ are pairwise disjoint.
We will build the final ratio equal to $\rho^+_i$ by performing a perturbation of $g$ with support in  $ g^{-w_i} (J_i),\dots,J_i$,  and  the initial ration equal to $\rho^-_i$ by a  perturbation of $g$ with support in  $I_i, g(I_i),\dots, g^{w_i}(I_i)$.  This supports are disjoint so that the construction can be performed independently.  Furthermore, they are analogous.  Will will present only the construction of $\rho^-_i$.
\vskip 2mm

For any $j$ the restriction $g|_{g^j(I_i)}$ maps $g^j(I_i)$ to $g^{j+1}(I_i)$. It will be more confortable to deal with diffeomorphisms of the same interval.  For that we will rescale the intervals $g^j(I_i)$ by affine maps on $[0,1]$. As this rescaling is affine it will not affect   the distorsion of $g|_{g^j(I_i)}$, and small $C^1$-pertubation of the rescaled map will induce $C^1$-pertubation of $g$ of proportional $C^1$-size.  More precisely:

\vskip 2mm

Let $\varphi_{i,j}\colon g^j(I_i(g))\to [0,1] $, $j\in\{0,\dots,w_i\}$ be the affine orientation preserving maps. We denote by $G_{i,j}\colon [0,1]\to [0,1]$, $j\in\{0,\dots,w_i-1\}$ the diffeomorphism $\varphi_{i,j+1} \circ g|_{g^j(I_i)}\circ \varphi_{i,j}^{-1}$.

As $g$ is $C^2$ the orbits are all dense so that the length of the interval ${g^j(I_i)}$ tends uniformly to $0$ when $i$ tends to infinity.  As a consequence the distorsion of $g$ tends to $0$ on ${g^j(I_i)}$.  As a direct consequence the diffeomorphisms $G_{i,j}$ tend uniformly to the identity map in the $C^1$ topology.

Notice that, according to Remark~\ref{r.subexp}, there are $\tilde\lambda_i>1$ tending to $1$ so that the initial and final ratio of $g$ belong to $[\tilde\lambda_i^{-w_i},\tilde\lambda_i^{w_i}]$.

Furthermore, as $g$ is assumed here to be $C^2$ and as the intervals $g^j(I_i)$ are pairwise disjoint in $S^1$, so that the sum of their length is bounded by $1$, a very classical distorsion control argument implies  that
\begin{lemm} There is a constant $C>1$ such that for every $i$, for any $ j\in \{0,1,..., w_i\}$ the distorsion of $g^{j}$ on $I_i$ is bounded by $\log C$.
\end{lemm}

As a consequence one gets that for every $i$ one has $$\frac{|g^{r(k_i)+j}(y)-g^j(y)|}{|g^{s(k_i)+j}(y)-g^j(y)|}\in [C^{-2}\tilde\lambda_i^{-w_i},C^2 \tilde\lambda_i^{w_i}]$$

Up to replace   $\lambda_i$ and  $\tilde \lambda_i$ by $\mu_i=\sup \{\lambda_i, C^{\frac 2{w_i}}\tilde\lambda_i\}$ one has that $ \mu_i \to 1$ and the numbers

$$\rho^-_i, \frac{|g^{r(k_i)+j}(y)-g^j(y)|}{|g^{s(k_i)+j}(y)-g^j(y)|}\in [\mu_i^{-w_i},\mu_i^{w_i}].$$

Our main lemma is
\begin{lemm}\label{l.perturbation} Let $w_i$ be a sequence tending to infinity. Let $G_{i,j}$, $j\in\{0, \dots,w_i-1\}$ be families of diffeomorphisms of $[0,1]$ tending uniformly to identity in the $C^1$ topology with $i$.

Then given any  points $t_{i,1}$, $t_{i,2}$ satisfying $\frac{t_{i,\eta}}{1-t_{i,\eta}}\in [\mu_i^{-w_i},\mu_i^{w_i}]$, $\eta\in\{1,2\}$, and given any $\varepsilon>0$. Then there is $i_0$ so that for any $i\geq i_0$, there are families $H_{i,j}$ so that:
\begin{itemize}
\item $H_{i,j}$ are $\varepsilon$ $C^1$-close to $G_{i,j}$
\item $H_{i,j}$ coincides with $G_{i,j}$ in neighborhoods of $0$ and $1$
\item $H_{i,w_i-1}\circ \cdots\circ H_{i,0}(t_{i,1})=t_{i,2}$.
\end{itemize}

\end{lemm}

\begin{demo}[Proof of Proposition~\ref{p.caracteristic}] Consider $g_1$ defined has $g_1=g$ out of the union of the intervals $g^j(I_i)$, $j\in\{0,\dots,w_i-1\}$ and equal to $\varphi_{i,j+1}^{-1} \circ H_{i,j}\circ \varphi_{i,j}$ on $g^j(I_i)$, where $H_{i,j}$ is given by Lemma~\ref{l.perturbation} for the constant:
\begin{itemize}
\item  $\frac \varepsilon M$ where $M$ is a bound for $dg$.
\item $t_2$ is $\varphi_{i,w_i}(g^{w_i}(y))$
\item $t_1$ is the point such that $\frac{t_1}{1-t_1}$ is the initial ration $\rho^-_i$.
\end{itemize}

Using the fact that $g=\varphi_{i,j+1}^{-1} \circ G_{i,j}\circ \varphi_{i,j}$ on $g^j(I_i)$ and that $H_{i,j}=G_{i,j}$ is a neighborhood of $0$ and $1$ one easily check that $g_1$ is a diffeomorphism.  Furthermore, the fact that $H_{i,j}$ is $\varepsilon/M$ close to $G_{i,j}$ implies that $g_1$ is $\varepsilon$ close to $g$ is the $C^1$ topology.  Finally the orbit segment of length $k_i$ through the point $y_1=\varphi_{i,0}^{-1}(t_1)$ satisfies
\begin{itemize}
\item   $g_1^{w_i}(y_1)=g^{w_i}(y)$ and therefore, for any $j\in\{w_i,\dots,k_i\}$  $g_1^j(y_1)=g^j(y)$
\item $\{y_1,\dots, g_1^{k_i}(y_1)\}$ is an adapted segment for $g_1$ ordered as the adapted segment $\{y,\dots, g^{k_i}(y)\}$
\end{itemize}

One deduces that the initial ratio of this segment is $\rho^-_i$, as announced.

\end{demo}

\subsection{Proof of Lemma~\ref{l.perturbation}}

Notice that, as the $G_{i,j}$ are assumed to tend uniformly to identity, the condition that the $H_{i,j}$ are $\varepsilon$-close to the $G_{i,j}$ can be replaced by the condition $H_{i,j}$ $\varepsilon$ $C^1$-close to identity (up to shrink slightly $\varepsilon$). Furthermore, the condition that $H_{i,j}$ and $G_{i,j}$ coincide in an (arbitrarily small) neighborhood of $0$ and $1$ can be obtained by the use of a bump function, without introducing derivatives larger that $1+2\varepsilon$.  Therefore, up to replace $\varepsilon$ by $\varepsilon/2$ Lemma~\ref{l.perturbation} is a direct consequence of the following:
\begin{lemm}\label{l.perturbation2} Let $w_i$ be a sequence tending to infinity and $\mu_i$ be a sequence tending to one.

Then given any  points $t_{i,1}$, $t_{i,2}$ satisfying $\frac{t_{i,\eta}}{1-t_{i,\eta}}\in [\mu_i^{-w_i},\mu_i^{w_i}]$, $\eta\in\{1,2\}$, and given any $\varepsilon>0$. Then there is $i_0$ so that for any $i\geq i_0$, there are families $H_{i,j}$ so that:
\begin{itemize}
\item $H_{i,j}$ are $\varepsilon$ $C^1$-close to $Id$
\item $H_{i,w_i-1}\circ \cdots\circ H_{i,0}(t_{i,1})=t_{i,2}$.
\end{itemize}
\end{lemm}

The main step for the proof is the elementary following

\begin{lemm}\label{l.1step} Given $\varepsilon>0$ small enough and $t\in[-1,1]$ and $y \in[0,1]$ there is a diffeomorphism $\varphi$ of $[0,1]$ which is equal to identity in a neighborhood of $0$ and $1$,  $2\varepsilon$-$C^1$-close to the identity, and such that

$$ \frac{|\varphi(y)|}{|\varphi(y)-\varphi(1)|}= (1+t\varepsilon) \frac{|y|}{|1-y|}$$
\end{lemm}
\begin{demo} Let $y_1$ be the point in $[0,1]$ so that $\frac{y_1}{1-y_1}=(1+t\varepsilon) \frac{|y|}{|1-y|}$.  An easy calculation shows that
$y_1= \frac{(1+t\varepsilon)y}{1+t\varepsilon y}$

The map $\varphi$ is just obtained by smoothing the piecewise affine homeomorphism, affine from $[0,y]$ to $[0,y_1]$ and from $[y,1]$ to $[y_1,1]$.  Notice that the linear parts of the affine segments are
\begin{itemize}
\item $\frac {y_1}y= \frac{1+t\varepsilon}{1+t\varepsilon y}=1+ t\varepsilon\frac{1-y}{1+t\varepsilon y}\in (1-\varepsilon, 1+\varepsilon)$  and
\item $\frac {1-y_1}{1-y}=\frac 1{1+t\varepsilon y}=1-\left(\frac{t y}{1+t\varepsilon y}\right)\varepsilon\in (1-2\varepsilon, 1 +2\varepsilon)$ for $\varepsilon <\frac{1}{2}$.
\end{itemize}

\end{demo}

\begin{demo}[Proof of Lemma~\ref{l.perturbation2}] Applying Lemma~\ref{l.1step}   $w_i$ times to $\frac{\varepsilon}{2}$, it is enough that
$\frac{t_{i,2}}{1-t_{i,2}}= (1+t\frac{\varepsilon}{2})^{w_i}\frac{t_{i,1}}{1-t_{i,1}}$ for some $t\in [0,1]$.

By assumption  $\frac{t_{i,1}}{1-t_{i,1}}, \frac{t_{i,2}}{1-t_{i,2}}\in [\mu_i^{-w_i},\mu_i^{w_i}]$, that is
$\frac{t_{i,2}}{1-t_{i,2}} \frac{1-t_{i,1}}{t_{i,1}}\in [\mu_i^{-2w_i},\mu_i^{2w_i}]$

Therefore, one can find $t$ if $\mu_i^2< 1+\frac{\varepsilon}{2}$. As $\mu_i\to 1$ when $i\to\infty$, it is enough to choose $i$ large enough, ending the proof.
\end{demo}